\documentclass{article}
\usepackage{amssymb,latexsym,amsmath,amsthm,mathrsfs}
\usepackage{graphicx}
\newcommand{\comment}[1]{}

\begin{document}
\title{On the sum of the series formed from the prime numbers
$\frac{1}{3}-\frac{1}{5}+
\frac{1}{7}+\frac{1}{11}-\frac{1}{13}-\frac{1}{17}+\frac{1}{19}
+\frac{1}{23}-\frac{1}{29}+\frac{1}{31}-$ etc.
where the prime numbers of the form $4n-1$ have a positive sign and those of the form $4n+1$ a negative sign\footnote{Presented to the
St. Petersburg Academy on October 2, 1775.
Originally published as
{\em De summa seriei ex numeris primis formatae $\frac{1}{3}-\frac{1}{5}+
\frac{1}{7}+\frac{1}{11}-\frac{1}{13}-\frac{1}{17}+\frac{1}{19}+\frac{1}{23}-\frac{1}{29}+\frac{1}{31}-$ etc. ubi numeri primi formae $4n-1$ habent signum positivum formae autem $4n+1$ signum negativum},
Opuscula analytica \textbf{2} (1785), 240--256.
E596 in the Enestr{\"o}m index.
Translated from the Latin by Jordan Bell,
Department of Mathematics, University of Toronto, Toronto, Ontario, Canada.
Email: jordan.bell@gmail.com}}
\author{Leonhard Euler}
\date{}
\maketitle

1. Even as Euclid had demonstrated that the multitude of prime numbers is infinite,
 many years ago I also showed that the sum of the series of the reciprocals
of the primes, namely
\[
\frac{1}{2}+\frac{1}{3}+\frac{1}{5}+\frac{1}{7}+\frac{1}{11}+\frac{1}{13}+\frac{1}{17}+\frac{1}{19}+\textrm{etc.}
\]
is infinitely large; more precisely, I showed that it has the magnitude of
the logarithm of the harmonic series
\[
1+\frac{1}{2}+\frac{1}{3}+\frac{1}{4}+\frac{1}{5}+\textrm{etc.},
\]
which seems not just a little remarkable, since commonly the harmonic
series is counted as the smallest kind of the infinite. However, since
not only is the logarithm of an infinite number itself infinite, but also
the logarithms of these logarithms are even still infinite,
and this clearly gives infinitely many lower degrees of the infinite.
Thus if $A$ denotes the sum of the series of reciprocals of the prime numbers,
then $lA$ will be infinitely large but it will also be considered to
belong to an infinitely lower order of the infinite; 
then indeed even now these formulae: $llA,lllA,llllA$ etc. will be infinite,
though each of these will be infinitely smaller than the preceding.

2. We next observe that all the prime numbers, aside from $2$, 
can be naturally divided into two classes, as those of the form $4n+1$ and those
of the form $4n-1$, where all the former are the sum of two squares while
the latter are entirely excluded from having this property.
The series formed by the reciprocals from these two classes are:
\[
\frac{1}{5}+\frac{1}{13}+\frac{1}{17}+\frac{1}{29}+\textrm{etc.} \quad
\textrm{and} \quad \frac{1}{3}+\frac{1}{7}+\frac{1}{11}+\frac{1}{19}+
\frac{1}{23}+\textrm{etc.}
\]
and even both are infinite, which also holds for all types of prime numbers.
Thus if from the prime numbers those are selected which are of the form
$100n+1$, of which sort are $101,401,601,701$ etc., not only will the multitude of them be infinite, but also
 the sum of this series formed from them, namely:
\[
\frac{1}{101}+\frac{1}{401}+\frac{1}{601}+\frac{1}{701}+\frac{1}{1201}+
\frac{1}{1301}+\frac{1}{1601}+\frac{1}{1801}+\frac{1}{1901}+\textrm{etc.},
\]
will also be infinite.

3. Let us first consider here the difference between the prime numbers 
of the forms $4n+1$ and $4n-1$, for there would be no doubt that 
the series formed from each of these 
are infinite should the difference between them have a definite value.
For this let us take the
terms of the form $4n-1$ to have the sign $+$,
and indeed the remaining terms to have the sign $-$,
so that this series will arise:
\[
\frac{1}{3}-\frac{1}{5}+\frac{1}{7}+\frac{1}{11}-\frac{1}{13}-\frac{1}{17}+\frac{1}{19}+\frac{1}{23}-\frac{1}{29}+\frac{1}{31}-\textrm{etc.}
\]
in which there one discerns no clear order for the signs. We shall see however
this does not hinder us from at least assigning an approximate sum. 
It does not seem unlikely that this sum, if it is not rational or irrational,
would in any case belong to a very notable kind of transcendental quantity.
In the meantime however, as not only is no clear order apparent in the signs, but even
less in the fractions themselves, at first sight no way seems open 
by which the sum of this can be reached.

4. In contemplating the well known Leibniz series for giving the quadrature of the circle
\[
\frac{\pi}{4}=1-\frac{1}{3}+\frac{1}{5}-\frac{1}{7}+\frac{1}{9}-\frac{1}{11}+\textrm{etc.},
\]
we see in it that all the odd numbers of the form $4n+1$ have the sign $+$,
and the remaining, of course of the form $4n-1$, have the sign $-$. All the
terms of our series with their corresponding
signs occur in the
Leibniz series with the signs changed. Then, if we remove all the composite numbers
from this series, finally of series, besides unity, would remain, with
opposite signs.
Thus we will obtain the sum of our series if we successively remove
all the composite numbers from the Leibniz series, since
the sum of all the
terms which
are excluded by one of these operations can be easily assigned.

5. Let us therefore begin with the Leibniz series itself, setting 
\[
A=1-\frac{1}{3}+\frac{1}{5}-\frac{1}{7}+\frac{1}{9}-\frac{1}{11}+\frac{1}{13}-\frac{1}{15}+\frac{1}{17}-\frac{1}{19}+\textrm{etc.},
\]
so $A=\frac{\pi}{4}$,
and let us first exclude all composite numbers divisible by $3$.
To this end, let us form this series:
\[
(A-1)\frac{1}{3}=-\frac{1}{9}+\frac{1}{15}-\frac{1}{21}+\frac{1}{27}-\frac{1}{33}+\frac{1}{39}-\textrm{etc.}
\]
which certainly contains all the composite numbers divisible by $3$, so if
this series is added to the former, all these composite terms will be excluded
\[
\frac{4}{3}A-\frac{1}{3}=1-\frac{1}{3}+\frac{1}{5}-\frac{1}{7}-\frac{1}{11}+\frac{1}{13}+\frac{1}{17}-\frac{1}{19}-\frac{1}{23}+\frac{1}{25}+\textrm{etc.}=B,
\]
so $B=\frac{4}{3}A-\frac{1}{3}$. Therefore in this series, whose sum
we have found, no more composite numbers divisible by $3$ occur, but
the first composite term occurring in it is $+\frac{1}{25}$.

6. Therefore let us now exclude all the terms divisible by $5$ from this last
series, to which end let us form this series:
\[
(B-1+\frac{1}{3})\frac{1}{5}=\frac{1}{25}-\frac{1}{35}-\frac{1}{55}+\frac{1}{65}+\frac{1}{85}-\frac{1}{95}-\textrm{etc.},
\]
which includes all terms divisible by five, 
which were so far still included in the series $B$, with the same
signs. With this latter series removed from the former, the following will remain:
\[
\frac{4}{5}B+\frac{1}{5}(1-\frac{1}{3})=1-\frac{1}{3}+\frac{1}{5}-\frac{1}{7}-\frac{1}{11}+\frac{1}{13}+\frac{1}{17}-\frac{1}{19}-\frac{1}{23}+\textrm{etc.},
\]
where now from the start all the terms are prime, and the first term
which occurs that is not a prime will be $\frac{1}{49}$, the next indeed
$\frac{1}{77},-\frac{1}{91}$.
Let us also put the sum of this series $=C$, so
\[
C=\frac{4}{5}B+\frac{1}{5}(1-\frac{1}{3}).
\]

7. Now let us therefore exclude all the terms which are still divisible by seven,
which the following form contains:
\[
(C-1+\frac{1}{3}-\frac{1}{5})\frac{1}{7}=-\frac{1}{49}-\frac{1}{77}+\frac{1}{91}+\textrm{etc.},
\]
in which the signs of the terms are opposite. 
This series added to the series $C$ gives
\[
\frac{8}{7}C-\frac{1}{7}(1-\frac{1}{3}+\frac{1}{5})=1-\frac{1}{3}+\frac{1}{5}-\frac{1}{7}-\frac{1}{11}+\frac{1}{13}+\frac{1}{17}-\textrm{etc.},
\]
in which the first composite term will be $+\frac{1}{121}$,
the next indeed $-\frac{1}{11\cdot 13}, -\frac{1}{11\cdot 17}$ etc. 
As well, let us call this series $D$, so that
\[
D=\frac{8}{7}C-\frac{1}{7}(1-\frac{1}{3}+\frac{1}{5}).
\]

8. Now from the series $D$ which has just been found, let us remove
the terms which are still divisible by $11$, which this form contains:
\[
(D-(1-\frac{1}{3}+\frac{1}{5}-\frac{1}{7}))\frac{1}{11}=-\frac{1}{121}+\frac{1}{143}+\frac{1}{187}-\textrm{etc.},
\]
which terms in the series $D$ have opposite signs;
whence if this series is added to the former these terms will be excluded,
and it will follow
\[
\frac{12}{11}D-\frac{1}{11}(1-\frac{1}{3}+\frac{1}{5}-\frac{1}{7})=1-\frac{1}{3}+\frac{1}{5}-\frac{1}{7}+\frac{1}{11}+\textrm{etc.},
\]
in which the first term which is not prime is $\frac{1}{169}$;
let us further designate this series by the letter $E$, so that
\[
E=\frac{12}{11}D-\frac{1}{11}(1-\frac{1}{3}+\frac{1}{5}-\frac{1}{7}).
\]

9. Therefore from this series let us exclude all the terms in it which are still
divisible by $13$, which this form contains:
\[
(E-1+\frac{1}{3}-\frac{1}{5}+\frac{1}{7}+\frac{1}{11})\frac{1}{13}=\frac{1}{169}+\frac{1}{221}-\textrm{etc.},
\]
and these terms have the same sign as in the series $E$ itself.
Thus this series should be subtracted from the former,
from which follows
\[
\frac{12}{13}E+\frac{1}{13}(1-\frac{1}{3}+\frac{1}{5}-\frac{1}{7}-\frac{1}{11})=1-\frac{1}{3}+\frac{1}{5}-\frac{1}{7}-\frac{1}{11}+\textrm{etc.},
\]
where the first term which is not prime is $\frac{1}{289}$. 
As well, let us designate this whole series by the letter $F$, so that
\[
F=\frac{12}{13}E+\frac{1}{13}(1-\frac{1}{3}+\frac{1}{5}-\frac{1}{7}-\frac{1}{11}).
\]

10. Now, since if we continue these operations further,
as long as we successively exclude the terms still divisible by $17$,
then indeed by $19$, by $23$ etc., finally 
such a series of prime numbers coming after unity will remain, which,
if we designate with the letter $Z$,
which should be thought of as continued infinitely, will be
\[
Z
=1-\frac{1}{3}+\frac{1}{5}-\frac{1}{7}-\frac{1}{11}+\frac{1}{13}+\frac{1}{17}-\frac{1}{19}-\textrm{etc.},
\]
consequently the sum of the series proposed in the title will be $1-Z$.
And it is also clear that these formulae continually come closer to this value:
\[
1-A, \quad 1-B, \quad 1-C, \quad 1-D, \quad 1-E, \quad 1-F \quad \textrm{etc.}
\]

11. It will become clear from the following formulae how the values
of these letters are successively determined from the preceding:
\begin{eqnarray*}
B&=&\frac{4}{3}A-\frac{1}{3}\cdot 1\\
C&=&\frac{4}{5}B+\frac{1}{5}(1-\frac{1}{3})\\
D&=&\frac{8}{7}C-\frac{1}{7}(1-\frac{1}{3}+\frac{1}{5})\\
E&=&\frac{12}{11}D-\frac{1}{11}(1-\frac{1}{3}+\frac{1}{5}-\frac{1}{7})\\
F&=&\frac{12}{13}E+\frac{1}{13}(1-\frac{1}{3}+\frac{1}{5}-\frac{1}{7}-\frac{1}{11})\\
G&=&\frac{16}{17}F+\frac{1}{17}(1-\frac{1}{3}+\frac{1}{5}-\frac{1}{7}-\frac{1}{11}+\frac{1}{13})\\
H&=&\frac{20}{19}G-\frac{1}{19}(1-\frac{1}{3}+\frac{1}{5}-\frac{1}{7}-\frac{1}{11}+\frac{1}{13}+\frac{1}{17})\\
I&=&\frac{24}{23}H-\frac{1}{23}(1-\frac{1}{3}+\frac{1}{5}-\frac{1}{7}-\frac{1}{11}+\frac{1}{13}+\frac{1}{17}-\frac{1}{19})\\
&&\textrm{etc.}
\end{eqnarray*}
One should note here that if the denominator is a prime of the form $4n+1$ then
the numerator of the first term will be one less, or $4n$, and also
indeed the next term should be added. On the other hand, if the denominator
is a prime of the form $4n-1$, then the numerator of the first term will be
one greater, or $4n$, and indeed the next term should be subtracted in this case.

12. So that now we can express all these values in numbers by decimal
fractions, first of all it is noted to be
\[
A=\frac{\pi}{4}=0,7853981634.
\]
Also for the other letters the following values are computed:
\[
\begin{array}{ll}
1-\frac{1}{3}&=b=0,6666666666\\
1-\frac{1}{3}+\frac{1}{5}&=c=0,8666666666\\
1-\frac{1}{3}+\frac{1}{5}-\frac{1}{7}&=d=0,7238095238\\
1-\frac{1}{3}+\frac{1}{5}-\frac{1}{7}-\frac{1}{11}&=e=0,6329004329\\
1-\frac{1}{3}+\frac{1}{5}-\frac{1}{7}-\frac{1}{11}+\frac{1}{13}&=f=0,7098235098\\
1-\frac{1}{3}+\frac{1}{5}-\frac{1}{7}-\frac{1}{11}+\frac{1}{13}+\frac{1}{17}&=g=0,7686470392\\
1-\frac{1}{3}+\frac{1}{5}-\frac{1}{7}-\frac{1}{11}+\frac{1}{13}+\frac{1}{17}-\frac{1}{19}&=h=0,7160154603\\
1-\frac{1}{3}+\frac{1}{5}-\frac{1}{7}-\frac{1}{11}+\frac{1}{13}+\frac{1}{17}-\frac{1}{19}-\frac{1}{23}&=i=0,6725371994,
\end{array}
\]
in which order the first term $a$ is equal to unity.

13. For the computation of these values $A,B,C,D,E$ etc. it will be
helpful to use the following formulae, in which at once we shall write the numerical
values of these letters
\[
\begin{array}{rclcl}
B&=&A+\frac{1}{3}(A-a)&=&0,713864\\
C&=&B-\frac{1}{5}(B-b)&=&0,704424\\
D&=&C+\frac{1}{7}(C-c)&=&0,681247\\
E&=&D+\frac{1}{11}(D-d)&=&0,677377\\
F&=&E-\frac{1}{13}(E-e)&=&0,673956\\
G&=&F-\frac{1}{17}(F-f)&=&0,676066\\
H&=&G+\frac{1}{19}(G-g)&=&0,671193\\
I&=&H+\frac{1}{23}(H-h)&=&0,699245\\
K&=&I-\frac{1}{29}(I-i)&=&0,669358.
\end{array}
\] 

14. As far as we have led this calculation, still we cannot be certain
beyond the third decimal figure of the sum of our series, and thus
it remains in doubt whether this sum is either greater or less than
$0,669$.
Now if we assume this value as true, 
the given series
\[
\frac{1}{3}-\frac{1}{5}+\frac{1}{7}+\frac{1}{11}-\frac{1}{13}-\frac{1}{17}+\frac{1}{19}+\textrm{etc.}
\]
will have the sum $0,331$,
and this value will be a bit smaller than $\frac{1}{3}$.
Since however by removing $\frac{1}{5}$ and then adding
$\frac{1}{7}+\frac{1}{11}$, as the sum of these fractions is greater
than $\frac{1}{5}$, it might happen that the true value
exceeds $\frac{1}{3}$, which at this point remains in doubt.
Indeed, for inquiring into the sum of this series another 
very accurate method is given, which we will expand here,
as it seems worth the effort 
to approximate the true sum of this series.

15. By the method in which we successively removed
the composite terms from the first Leibniz series, likewise if
we remove entirely all the terms besides unity we will
find
\[
\frac{\pi}{4}=\frac{3}{4}\cdot \frac{5}{4}\cdot \frac{7}{8}\cdot \frac{11}{12}
\cdot \frac{13}{12} \cdot \frac{17}{16}\cdot \frac{19}{20} \cdot \textrm{etc.},
\]
where in the numerators all the prime numbers except $2$ occur,
and indeed the denominators are the evenly even numbers, either greater than
or less than their numerator. Also indeed, if the reciprocal series of
odd squares:
\[
1+\frac{1}{3^2}+\frac{1}{5^2}+\frac{1}{7^2}+\frac{1}{9^2}+\frac{1}{11^2}
+\frac{1}{13^2}+\textrm{etc.},
\]
whose sum I have shown to be $=\frac{\pi\pi}{8}$,
is treated in a similar way, it will turn into
\[
\frac{\pi\pi}{8}=\frac{3\cdot 3}{2\cdot 4}\cdot \frac{5\cdot 5}{4\cdot 6}
\cdot \frac{7\cdot 7}{6\cdot 8}\cdot \frac{11\cdot 11}{10\cdot 12}\cdot \frac{13\cdot 13}{12\cdot 14}\cdot \textrm{etc.},
\]
where again in the numerators all the prime numbers occur twice,
and in the denominators the same again, either increased or decreased by one.
Thus if we divide this expression by the square of the former, which is
\[
\frac{\pi\pi}{16}=\frac{3\cdot 3}{4\cdot 4}\cdot \frac{5\cdot 5}{4\cdot 4}\cdot
\frac{7\cdot 7}{8\cdot 8}\cdot \frac{11\cdot 11}{12\cdot 12}\cdot
\frac{13\cdot 13}{12\cdot 12}\cdot \textrm{etc.},
\]
the quotient will be
\[
2=\frac{4}{2}\cdot \frac{4}{6}\cdot \frac{8}{6}\cdot \frac{12}{10}
\cdot \frac{12}{14}\cdot \textrm{etc.},
\]
where all the prime numbers occur either increased or decreased by unity,
and where the evenly even numbers appear in the numerator and the oddly even
in the the denominator.

16. Next this expression can then be exhibited in this way:
\[
2=\frac{3+1}{3-1}\cdot \frac{5-1}{5+1}\cdot \frac{7+1}{7-1}\cdot \frac{11+1}{11-1}\cdot \frac{13-1}{13+1}\cdot \textrm{etc.};
\]
then by taking the hyperbolic logarithm we will have:
\[
l2=l\frac{3+1}{3-1}+l\frac{5-1}{5+1}+l\frac{7+1}{7-1}+l\frac{11+1}{11-1}+l\frac{13-1}{13+1}+\textrm{etc.}
\]
It is known however in general by infinite series that
\[
\frac{1}{2}l\frac{a+1}{a-1}=\frac{1}{a}+\frac{1}{3a^3}+\frac{1}{5a^5}+\frac{1}{7a^7}
+\frac{1}{9a^9}+\textrm{etc.},
\]
and so 
\[
\frac{1}{2}l\frac{a-1}{a+1}=-\frac{1}{a}-\frac{1}{3a^3}-\frac{1}{5a^5}-\frac{1}{7a^7}
-\frac{1}{9a^9}-\textrm{etc.}
\]
Now if by means of these formulae we convert all the former logarithms into
infinite series, even though innumerable infinite series will be obtained,
one will be able to reduce these to series which can be dealt with easily.

17. Thus, first all these logarithms should be divided by two,
and as the hyperbolic logarithm is being taken here, as
\[
l2=0,6931471805
\]
it will be
\[
\frac{1}{2}l2=0,3465735902,
\]
while from the other side the logarithms will be arranged thus:
\begin{eqnarray*}
\frac{1}{2}l\frac{3+1}{3-1}&=&\frac{1}{3}+\frac{1}{3\cdot 3^3}+\frac{1}{5\cdot 3^5}
+\frac{1}{7\cdot 3^7}+\frac{1}{9\cdot 3^9}+\textrm{etc.}\\
\frac{1}{2}l\frac{5-1}{5+1}&=&-\frac{1}{5}-\frac{1}{3\cdot 5^3}-\frac{1}{5\cdot 5^5}
-\frac{1}{7\cdot 5^7}-\frac{1}{9\cdot 5^9}-\textrm{etc.}\\
\frac{1}{2}l\frac{7+1}{7-1}&=&\frac{1}{7}+\frac{1}{3\cdot 7^3}
+\frac{1}{5\cdot 7^7}+\frac{1}{7\cdot 7^7}+\frac{1}{9\cdot 7^9}+\textrm{etc.}\\
\frac{1}{2}l\frac{11+1}{11-1}&=&\frac{1}{11}+\frac{1}{3\cdot 11^3}+\frac{1}{5\cdot 11^5}+\frac{1}{7\cdot 11^7}+\frac{1}{9\cdot 11^9}+\textrm{etc.}\\
\frac{1}{2}l\frac{13-1}{13+1}&=&-\frac{1}{13}-\frac{1}{3\cdot 13^3}-\frac{1}{5\cdot 13^5}-\frac{1}{7\cdot 13^7}-\frac{1}{9\cdot 13^9}-\textrm{etc.}\\
&&\textrm{etc.}
\end{eqnarray*}

18. Now by descending vertically, let us consider the following series which are
also infinite:
\begin{eqnarray*}
O&=&\frac{1}{3}-\frac{1}{5}+\frac{1}{7}+\frac{1}{11}-\frac{1}{13}-\frac{1}{17}+\frac{1}{19}+\textrm{etc.}\\
P&=&\frac{1}{3^3}-\frac{1}{5^3}-\frac{1}{7^3}+\frac{1}{11^3}-\frac{1}{13^3}-\frac{1}{17^3}+\frac{1}{19^3}+\textrm{etc.}\\
Q&=&\frac{1}{3^5}-\frac{1}{5^5}+\frac{1}{7^5}+\frac{1}{11^5}-\frac{1}{13^5}
-\frac{1}{17^5}+\frac{1}{19^5}+\textrm{etc.}\\
R&=&\frac{1}{3^7}-\frac{1}{5^7}+\frac{1}{7^7}+\frac{1}{11^7}-\frac{1}{13^7}-\frac{1}{17^7}+\frac{1}{19^7}+\textrm{etc.}\\
S&=&\frac{1}{3^9}-\frac{1}{5^9}+\frac{1}{7^9}+\frac{1}{11^9}-\frac{1}{13^9}-\frac{1}{17^9}+\frac{1}{19^9}+\textrm{etc.}\\
&&\textrm{etc.}
\end{eqnarray*}
The first of these series $O$ is the same one whose sum has been
given to us to investigate.

19. Therefore with these series so designated by capital letters, we will have
this equation:
\[
\frac{1}{2}l2=O+\frac{1}{3}P+\frac{1}{5}Q+\frac{1}{7}R+\frac{1}{9}S+\frac{1}{11}T+\textrm{etc.},
\]
whence if the sums of the series $P,Q,R,S$ etc. were known, from
them we would easily obtain the sum of the series $O$ that is being sought;
for it would be
\[
O=\frac{1}{2}l2-\frac{1}{3}P-\frac{1}{5}Q-\frac{1}{7}R-\frac{1}{9}S-\textrm{etc.}
\]

20. And we can conclude the sums of the series $P,Q,R$ etc. from the
class of series where all the odd numbers occur, in the same way in which
above we elicited the sum $O$ from the Leibniz series
\[
1-\frac{1}{3}+\frac{1}{5}-\frac{1}{7}+\frac{1}{9}-\frac{1}{11}+\textrm{etc.}
\]
To this end it will be useful to expand the following class of series:
\begin{eqnarray*}
\mathfrak{P}&=&1-\frac{1}{3^3}+\frac{1}{5^3}-\frac{1}{7^3}+\frac{1}{9^3}-\frac{1}{11^3}+\frac{1}{13^3}-\frac{1}{15^3}+\textrm{etc.}\\
\mathfrak{Q}&=&1-\frac{1}{3^5}+\frac{1}{5^5}-\frac{1}{7^5}+\frac{1}{9^5}-\frac{1}{11^5}+\frac{1}{13^5}-\frac{1}{15^5}+\textrm{etc.}\\
\mathfrak{R}&=&1-\frac{1}{3^7}+\frac{1}{5^7}-\frac{1}{7^7}+\frac{1}{9^7}-\frac{1}{11^7}+\frac{1}{13^7}-\frac{1}{15^7}+\textrm{etc.}\\
\mathfrak{S}&=&1-\frac{1}{3^9}+\frac{1}{5^9}-\frac{1}{7^9}+\frac{1}{9^9}-\frac{1}{11^9}+\frac{1}{13^9}-\frac{1}{15^9}+\textrm{etc.}\\
\mathfrak{T}&=&1-\frac{1}{3^{11}}+\frac{1}{5^{11}}-\frac{1}{7^{11}}+\frac{1}{9^{11}}-\frac{1}{11^{11}}+\frac{1}{13^{11}}-\frac{1}{15^{11}}+\textrm{etc.}\\
&&\textrm{etc.}
\end{eqnarray*}
But I have given the sums of all these series some time ago by the quadrature
of the circle, namely expressed by similar powers of $\pi$,
in the following way:
\[
\begin{array}{lll}
\mathfrak{P}=\frac{1}{1\cdot 2}\cdot \frac{\pi^3}{2^4}&&\mathfrak{T}=\frac{50521}{1\ldots 10} \frac{\pi^{11}}{2^{12}}\\
\mathfrak{Q}=\frac{5}{1\cdot 2\cdot 3\cdot 4}\cdot \frac{\pi^5}{2^6}&&\mathfrak{U}=\frac{2702765}{1\ldots 12}\cdot \frac{\pi^{13}}{2^{14}}\\
\mathfrak{R}=\frac{61}{1\ldots 6}\cdot \frac{\pi^7}{2^8}&&\mathfrak{V}=\frac{199360981}{1\ldots 14}\cdot \frac{\pi^{15}}{2^{16}}\\
\mathfrak{S}=\frac{1385}{1\ldots 8}\cdot \frac{\pi^9}{2^{10}}&&
\mathfrak{W}=\frac{19391512145}{1\ldots 16}\cdot \frac{\pi^{17}}{2^{18}}\\
&\textrm{etc.}&
\end{array}
\]

21. Let us thus expand these values into decimal fractions to six figures,
and it will be
\[
\begin{array}{lcl||l}
&&&\textrm{Differences}\\
\mathfrak{P}&=&0,9689462&0,0272116\\
\mathfrak{Q}&=&0,9961578&0,0033969\\
\mathfrak{R}&=&0,9995547&0,0003952\\
\mathfrak{S}&=&0,9999499&0,0000448\\
\mathfrak{T}&=&0,9999947&0,0000050\\
\mathfrak{U}&=&0,9999997&0,0000005\\
&&\textrm{etc.}&\textrm{etc.}
\end{array}
\] 

22. To now get at the values of the letters $P,Q,R$ etc.,
the same method can be used by which above we dismissed all the composite
terms from the series
\[
1-\frac{1}{3}+\frac{1}{5}-\frac{1}{7}+\frac{1}{9}-\frac{1}{11}+\textrm{etc.},
\]
except that in place of all the single numbers, powers of them should be written.
We shall explain these operations for the letters in general.
Thus let us consider this series:
\[
\mathfrak{Z}=1-\frac{1}{3^n}+\frac{1}{5^n}-\frac{1}{7^n}+\frac{1}{9^n}-
\frac{1}{11^n}+\textrm{etc.},
\]
whose sum, as above, we shall designate with the letter $A$, so that
$A=\mathfrak{Z}$, and then let us elicit the following letters $B,C,D$
etc. by the following formulae:
\[
\begin{array}{lll}
B=A+\frac{1}{3^n}(A-a)&\textrm{taking}&a=1\\
C=B-\frac{1}{5^n}(B-b)&&b=1-\frac{1}{3^n}\\
D=C+\frac{1}{7^n}(C-c)&&c=1-\frac{1}{3^n}+\frac{1}{5^n}\\
E=D+\frac{1}{11^n}(D-d)&&d=1-\frac{1}{3^n}+\frac{1}{5^n}-\frac{1}{7^n}\\
F=E-\frac{1}{13^n}(E-e)&&e=1-\frac{1}{3^n}+\frac{1}{5^n}-\frac{1}{7^n}-\frac{1}{11^n}\\
\textrm{etc.}&&\textrm{etc.}
\end{array}
\]
Having found these values, their complements to unity, namely
$1-A,1-B,1-C,1-D$ etc., will very quickly approach the sought for value
\[
Z=\frac{1}{3^n}-\frac{1}{5^n}+\frac{1}{7^n}+\frac{1}{11^n}-\frac{1}{13^n}-\textrm{etc.}
\]

23. Thus let us first apply these general precepts to the value of the letter
$P$,
for which we will start from the value
\[
\mathfrak{P}=0,9689462=A,
\] 
and because here it is $n=3$, we will have
\[
a=1, \quad b=0,9629630, \quad c=0,9709630, \quad d=0,9680476;
\]
additional values will not be useful. We will thus find the following values:
\[
\begin{array}{lclcl}
B&=&A-\frac{1}{3^3}\cdot 0,0310538&=&0,9677961\\
C&=&B-\frac{1}{5^3}\cdot 0,0048331&=&0,9677574\\
D&=&C-\frac{1}{7^3}\cdot 0,0032056&=&0,9677481\\
E&=&D-\frac{1}{11^3}\cdot 0,0002995&=&0,9677479.
\end{array}
\]
It will not be useful to go further; we thus will have 
\[
P=1-E=0,0322521,
\]
whence we now find that
\[
\frac{1}{2}l2-\frac{1}{3}\cdot 0,0322521=0,3358229.
\]

24. Let us now take $n=5$, and we will have
\[
A=\mathfrak{Q}=0,9961578,
\]
then indeed it will be
\[
a=1, \quad b=0,9958847, \quad c=0,9962048, \quad d=0,9961453,
\]
from which we will thus find
\[
\begin{array}{lclcl}
B&=&A-\frac{1}{3^5}\cdot 0,0038422&=&0,9961420\\
C&=&B-\frac{1}{5^5}\cdot 0,0002573&=&0,9961419.
\end{array}
\]
Thus it will be
\[
Q=1-C=0,0038581
\]
and hence
\[
\frac{1}{2}l2-\frac{1}{3}P-\frac{1}{5}Q=0,3350513.
\]

25. Now let $n=7$, and $A=\mathfrak{R}=0,9995547$, then indeed
$a=1,b=0,9995428$,  
and it will therefore become
\[
B=A-\frac{1}{3^7}\cdot 0,0004453=0,9995545,
\]
from which we will now have
\[
R=1-B=0,0004455
\]
and hence
\[
\frac{1}{2}l2-\frac{1}{3}P-\frac{1}{5}Q-\frac{1}{7}R=0,3349877.
\]

26. Since in this calculation it will nearly be $B=A$, in the following
the letter $B$ will not be useful, and therefore we will have
\[
S=1-\mathfrak{S}=0,0000501
\]
and so
\[
\frac{1}{9}S=0,0000056.
\]
Then indeed it will be
\[
T=1-\mathfrak{T}=0,0000053
\]
and so
\[
\frac{1}{11}T=0,0000005,
\]
and finally 
\[
U=1-\mathfrak{U}=0,0000003 \quad \textrm{and} \quad \frac{1}{13}U=0,0000000.
\]
Thus removing these parts from the preceding value gives
\[
O=0,3349816.
\]
It is apparent from this that the value is still somewhat greater than
$\frac{1}{3}$.

27. Thus we can now be certain that the sum of this infinite series
\[
\frac{1}{3}-\frac{1}{5}+\frac{1}{7}+\frac{1}{11}-\frac{1}{13}-\frac{1}{17}+\frac{1}{19}+\textrm{etc.}
\]
is almost exactly $=0,3349816$. It is now natural to investigate whether
or not this value holds some notable ratio either to the periphery of the circle
$\pi$ or to the hyperbolic logarithm of it, like we have observed above
that the reciprocal series of the prime numbers
\[
\frac{1}{2}+\frac{1}{3}+\frac{1}{5}+\frac{1}{7}+\frac{1}{11}+\textrm{etc.}
\]
expresses the hyperbolic logarithm of the complete harmonic series
\[
1+\frac{1}{2}+\frac{1}{3}+\frac{1}{4}+\frac{1}{5}+\frac{1}{6}+\textrm{etc.}
\]
It can likewise be seen that this series of prime numbers
\[
\frac{1}{3}-\frac{1}{5}+\frac{1}{7}+\frac{1}{11}-\frac{1}{13}-\frac{1}{17}+\textrm{etc.}
\]
also contains the logarithm of the compete series
\[
1-\frac{1}{3}+\frac{1}{5}-\frac{1}{7}+\frac{1}{9}-\frac{1}{11}+\textrm{etc.},
\]
whose sum is $\frac{\pi}{4}$. 
To this end I shall use the hyperbolic logarithm of $\pi$ which I have found
previously
\[
1,14472,98858,49400,17414,34273,51353,05865.
\]
Therefore we would like to see whether perhaps the sum turns out to be
$O=l\pi-lN$, where $N$ is a reasonably simple number. Truly though for the most part
the
investigations of this are undertaken without any 
success.

28. By means of the last method not only can we elicit the sum of
the given series, but also odd powers of it, and here we state
these sums for view:
\begin{eqnarray*}
\frac{1}{3}-\frac{1}{5}+\frac{1}{7}+\frac{1}{11}-\frac{1}{13}-\frac{1}{17}+\frac{1}{19}+\textrm{etc.}&=&0,3349816\\
\frac{1}{3^3}-\frac{1}{5^3}+\frac{1}{7^3}+\frac{1}{11^3}-\frac{1}{13^3}
-\frac{1}{17^3}+\frac{1}{19^3}+\textrm{etc.}&=&0,0322521\\
\frac{1}{3^5}-\frac{1}{5^5}+\frac{1}{7^5}+\frac{1}{11^5}-\frac{1}{13^5}-\frac{1}{17^5}+\frac{1}{19^5}+\textrm{etc.}&=&0,0038602\\
\frac{1}{3^7}-\frac{1}{5^7}+\frac{1}{7^7}+\frac{1}{11^7}-\frac{1}{13^7}-\frac{1}{17^7}+\frac{1}{19^7}+\textrm{etc.}&=&0,0004455\\
\frac{1}{3^9}-\frac{1}{5^9}+\frac{1}{7^9}+\frac{1}{11^9}-\frac{1}{13^9}-\frac{1}{17^9}+\frac{1}{19^9}+\textrm{etc.}&=&0,0000501\\
\frac{1}{3^{11}}-\frac{1}{5^{11}}+\frac{1}{7^{11}}+\frac{1}{11^{11}}-\frac{1}{13^{11}}-\frac{1}{17^{11}}+\frac{1}{19^{11}}+\textrm{etc.}&=&0,0000053\\
\frac{1}{3^{13}}-\frac{1}{5^{13}}+\frac{1}{7^{13}}+\frac{1}{11^{13}}-\frac{1}{13^{13}}-\frac{1}{17^{13}}+\frac{1}{19^{13}}+\textrm{etc.}&=&0,0000003\\
\textrm{etc.}&&
\end{eqnarray*}

29. It is doubtless that these sums deserve only a little attention unless
perhaps they can be reduced to known quantities. In truth, because 
in these series neither do the terms proceed according to a certain law,
or even in the signs is a more or less certain order seen, this inquiry
a first sight seems 
altogether impossible. Therefore the method by which we have obtained
these sums certainly is worth all attention, even more so
because it rests on quite intricate properties of series of powers.
For unless the sums of the series
\[
1-\frac{1}{3^n}+\frac{1}{5^n}-\frac{1}{7^n}+\frac{1}{9^n}-\textrm{etc.}
\]
in the cases when $n$ is an odd number are known, this whole investigation
may have been undertaken in vain.

\newpage

\section*{Translator's mathematical appendix}
\small{Define $\chi:\mathbb{Z} \to \{-1,0,1\}$ by
\[
\chi(n) = \begin{cases}
0&n \equiv 0 \pmod{4}\\
1&n \equiv 1 \pmod{4}\\
0&n \equiv 2 \pmod{4}\\
-1&n \equiv 3 \pmod{4}.
\end{cases}
\]
One checks that for any integers $m$ and $n$, $\chi(mn)=\chi(m) \chi(n)$, namely
$\chi$ is \textbf{completely multiplicative}. $\chi$ is the nontrivial Dirichlet character modulo 4.

We prove in this appendix that the series $\sum_p \frac{\chi(p)}{p}$ converges.
We follow an argument Davenport attributes to Mertens.\footnote{Harold Davenport,
{\em Multiplicative Number Theory}, third ed., p.~57, Chapter 7.}

We shall use \textbf{summation by parts}:
For $G_n=\sum_{k=1}^n g_k$, with $G_0=0$,
\begin{equation}
\sum_{k=m}^n f_k g_k = f_n G_n - f_m G_{m-1} - \sum_{k=m}^{n-1} G_k(f_{k+1}-f_k).
\label{parts}
\end{equation}

Define
\[
S(x) = \sum_{n \leq x} \chi(n).
\]
$S(x)$ satisfies
\[
S(1)=1, \quad S(2) = 1, \quad S(3) = 0, \quad S(4) = 0, \quad S(x+4)=S(x).
\]
Let $s>0$. 
We use summation by parts with
$f_n=n^{-s}$ and $g_n=\chi(n)$, so
$G_n=S(n)$. We get
\[
\sum_{n=1}^N \chi(n) n^{-s}  = N^{-s} S(N) - \sum_{n=1}^{N-1} S(n)((n+1)^{-s}-n^{-s}).
\]
Using the identity
\[
\int_n^{n+1} x^{-s-1} dx = \frac{x^{-s}}{-s} \big|_n^{n+1}= \frac{1}{-s}( (n+1)^{-s} - n^{-s}),
\]
and using that if $n \leq x < n+1$ then $S(n)=S(x)$, we get
\begin{align*}
\sum_{n=1}^N \chi(n) n^{-s}&=
N^{-s} S(N) +s \sum_{n=1}^{N-1} S(n)  \int_n^{n+1} x^{-s-1} dx\\
&=N^{-s} S(N) +s \sum_{n=1}^{N-1} \int_n^{n+1} S(x) x^{-s-1} dx\\
&=N^{-s} S(N) + s \int_1^N S(x) x^{-s-1} dx.
\end{align*}

For $s=1$,
\begin{align*}
\sum_{n=1}^N \chi(n) n^{-1}-\int_1^\infty S(x) x^{-2} dx&=
N^{-s} S(N) +  \int_1^N S(x) x^{-2} dx\\
&-\int_1^\infty S(x) x^{-2} dx\\
&=N^{-s} S(N)  -  \int_N^\infty S(x) x^{-2} dx.
\end{align*}
But
\[
\int_N^\infty S(x) x^{-2} dx\leq 
\int_N^\infty x^{-2} dx
=-x^{-1} \bigg|_N^\infty
=N^{-1}.
\]
Therefore,
\begin{equation}
\sum_{n=1}^N \chi(n) n^{-1}= A+ O(N^{-1}),
\qquad A=\int_1^\infty S(x) x^{-2} dx.
\label{errorterm}
\end{equation}

For a positive integer $n$ define
\[
\Lambda(n) = \begin{cases}
\log p&\textrm{$n$ is a power of a prime $p$}\\
0&\textrm{otherwise}.
\end{cases}
\]
For example, $\Lambda(1)=0, \Lambda(2)=\log 2, \Lambda(10)=0, \Lambda(27) = \log 3$.
$\Lambda$ is called the \textbf{von Mangoldt function}. One proves that it satisfies\footnote{G. H. Hardy
and E. M. Wright, {\em An Introduction to the Theory of Numbers}, fifth ed., p.~254, Theorem 296.}
\[
\log n = \sum_{m \mid n} \Lambda(m).
\]
For example,
\[
\sum_{m \mid 10} \Lambda(m) = \Lambda(10)+\Lambda(5)+\Lambda(2)+\Lambda(1)
=0+\log 5 + \log 2 + 0 = \log 10.
\]

Using the identity $\log n = \sum_{m \mid n} \Lambda(m)$, using that
$\chi$ is completely multiplicative, and using \eqref{errorterm},
\begin{align*}
\sum_{n \leq x} \frac{\chi(n) \log n}{n}&=\sum_{n \leq x, m_1 \mid n} \frac{\chi(n)}{n}  \Lambda(m_1)\\
&=\sum_{m_1 m_2 \leq x} \frac{\chi(m_1 m_2)}{m_1m_2} \Lambda(m_1)\\
&=\sum_{m_1 \leq x} \frac{\chi(m_1) \Lambda(m_1)}{m_1} \sum_{m_2 \leq \frac{x}{m_1}} \frac{\chi(m_2)}{m_2}\\
&=\sum_{m_1 \leq x} \frac{\chi(m_1) \Lambda(m_1)}{m_1} \left( A + O\big(\big(\frac{x}{m_1}\big)^{-1}\big)\right)\\
&=A \sum_{m_1 \leq x} \frac{\chi(m_1) \Lambda(m_1)}{m_1} + O\left(\frac{1}{x} \sum_{m_1 \leq x} \chi(m_1) \Lambda(m_1) \right).
\end{align*}
We have
\[
\left| \sum_{m_1 \leq x} \chi(m_1) \Lambda(m_1)\right| \leq \sum_{m_1 \leq x} |\chi(m_1) \Lambda(m_1)|
\leq \sum_{m_1 \leq x} \Lambda(m_1), 
\]
Let $\psi(x) = \sum_{m_1 \leq x} \Lambda(m_1)$. It satisfies the estimate
\[
\psi(x) = O(x),
\]
 due to Chebyshev.\footnote{G. H. Hardy
and E. M. Wright, {\em An Introduction to the Theory of Numbers}, fifth ed., p.~341, Theorem 414.}
Therefore,
\begin{align*}
\sum_{n \leq x} \frac{\chi(n) \log n}{n}&=A \sum_{m_1 \leq x} \frac{\chi(m_1) \Lambda(m_1)}{m_1} + O(1).
\end{align*}

\textbf{Dirichlet's test} states that if $a_n$ is a nonincreasing sequence of real numbers that tend to $0$ and 
$b_n$ is a sequence of complex numbers satisfying $\left| \sum_{n=1}^N b_n \right| \leq M$ for all $N$, then
the series $\sum_{n=1}^\infty a_n b_n$ converges. Applying Dirichlet's test, the series
$\sum_{n=1}^\infty \frac{\chi(n) \log n}{n}$ converges, and a fortiori, 
\[
\sum_{n \leq x}  \frac{\chi(n) \log n}{n} = O(1).
\]
It follows that
\[
A \sum_{m_1 \leq x} \frac{\chi(m_1) \Lambda(m_1)}{m_1} = O(1),
\]
and as $A>0$,
\begin{equation}
\sum_{n \leq x} \frac{\chi(n) \Lambda(n)}{n} = O(1).
\label{Oterm2}
\end{equation}

Let $g_n = \frac{\chi(n) \Lambda(n)}{n}$ and $f_n = \frac{1}{\log n}$. 
Applying summation by parts \eqref{parts},
\begin{align*}
\sum_{n=2}^N  \frac{\chi(n) \Lambda(n)}{n \log n}&=\frac{G_N}{\log N} - \frac{G_1}{\log 2}
-\sum_{n=2}^{N-1} G_n \left(\frac{1}{\log(n+1)}-\frac{1}{\log n}\right)\\
&=\frac{G_N}{\log N} + \sum_{n=2}^{N-1} G_n \left(\frac{1}{\log n}-\frac{1}{\log(n+1)}\right).
\end{align*}

Since $(\frac{1}{\log x})'=-\frac{1}{x(\log x)^2}$,
\[
\frac{1}{\log n}-\frac{1}{\log(n+1)} = \int_n^{n+1} \frac{1}{x(\log x)^2} dx,
\]
so, as $G_n=O(1)$,
\[
\left| G_n \left(\frac{1}{\log n}-\frac{1}{\log(n+1)}\right) \right|
=O\left(\frac{1}{n(\log n)^2}\right). 
\]
The series $\sum_{n=2}^\infty \frac{1}{n(\log n)^2}$ converges, so   
the series $\sum_{n=2}^\infty G_n \left(\frac{1}{\log n}-\frac{1}{\log(n+1)}\right)$ converges absolutely, and 
in particular it converges. It follows that the series $\sum_{n=2}^N  \frac{\chi(n) \Lambda(n)}{n \log n}$ converges.

We have
\begin{align*}
\sum_{n=2}^\infty  \frac{\chi(n) \Lambda(n)}{n \log n}&=\sum_p \sum_{k \geq 1}
\frac{\chi(p^k) \Lambda(p^k)}{p^k \log(p^k)}\\
&=\sum_p \frac{\chi(p)}{p} + \sum_p \sum_{k \geq 2} \frac{\chi(p^k) \Lambda(p^k)}{p^k  \log(p^k)}.
\end{align*}
But
\[
\left| \sum_{k \geq 2} \frac{\chi(p^k) \Lambda(p^k)}{p^k  \log(p^k)} \right|
\leq
\sum_{k \geq 2} p^{-k}
=\frac{1}{p(p-1)}.
\]
The series $\sum_p \frac{1}{p(p-1)}$ converges, and we have already established that the series $\sum_{n=2}^\infty  \frac{\chi(n) \Lambda(n)}{n \log n}$ converges. Therefore, the series
$\sum_p \frac{\chi(p)}{p}$ converges.}

\end{document}